\documentclass[11pt,a4paper]{amsart}
\usepackage[latin1]{inputenc}
\usepackage{amssymb}
\usepackage{amsmath}
\usepackage{latexsym}
\usepackage{xypic}
\usepackage{epic}
\usepackage{eepic}
\usepackage{graphicx}
\usepackage{psfrag}
\usepackage{longtable}
\usepackage{epsfig}
\usepackage{xypic}
\newtheorem{rem}{Remark}[section]
\newtheorem{prop}{Proposition}[section]
\newtheorem{lemma}{Lemma}[section]

\newtheorem{defi}{Definition}[section]
\newtheorem{theorem}{Theorem}[section]

\newcommand{\bprf}{{\it Proof.~}}
\newcommand{\eprf}{\hfill $\square$ \bigskip\par}

\newcommand{\ci}{ \mathbb{C}}
\newcommand{\R}{ \mathbb{R}}
\newcommand{\Z}{\mathbb{Z}}
\newcommand{\Q}{\mathbb{Q}}

\author{Alessandra~Sarti}

\begin{document}

\title{Transcendental lattices of some K3-Surfaces}

\address{\hskip -.43cm Alessandra Sarti, Fachbereich f\"ur Mathematik, Johannes Gutenberg-Universit\"at, 55099 Mainz, Germany}
\address{\hskip -.43cm {\it Current address}: Universit\`a di Milano, Dipartimento di Matematica, Via C. Saldini, 50, 20133 Milano, Italy}
\address{\hskip -.43cm {\it E-mail address}: {\tt sarti@mat.unimi.it}, {\tt sarti@mathematik.uni-mainz.de}}



\maketitle
\begin{abstract}
In a previous paper, \cite{io}, we described six families of $K3$-surfaces with Picard-number $19$, and we identified surfaces with Picard-number $20$. In these notes we classify some of the surfaces  by computing their transcendental lattices. Moreover we show that the surfaces with Picard-number $19$ are birational to a Kummer surface which is the quotient of a non-product type abelian surface by an involution.
\end{abstract}


\section{Introduction}\label{intro}
Given a K3-surface an important step toward its classification in view of the Torelli theorem is to compute the Picard lattice and the transcendental lattice. When the rank of the Picard lattice (i.e. the {\it Picard-number}, which we denote by $\rho$) of the K3-surface is 20, the maximal possible, the transcendental lattice has rank two. These $K3$-surfaces are called by Shioda and Inose {\it singular}. In  \cite{si}, Shioda and Inose classified such surfaces in terms of their transcendental lattice, more precisely they show the following:
\begin{theorem}\cite[Theorem 4, \S 4]{si}
There is a natural one-to-one correspondence from the set of singular $K3$-surfaces to the set of equivalence classes of positive-definite even integral binary quadratic forms with respect to $SL_2(\mathbb{Z})$.
\end{theorem}
When the Picard-number is $19$ the transcendental lattice has rank three and by results of Morrison, \cite{Mo}, and Nikulin, \cite{nikulin}, the embedding in the {\it K3-lattice} $\Lambda:=-E_8\oplus-E_8\oplus U\oplus U\oplus U$ is unique, hence it identifies the moduli curve classifying the K3-surfaces. In general however it seems to be difficult to compute explicitly the transcendental lattice. In \cite{io} we describe six families of K3-surfaces with Picard-number 19 and we identify in each family four surfaces with Picard-number 20. The aim of these notes is to compute their transcendental lattice and to classify them. In \cite{io} we describe completely the Picard lattice of the general surface in two of the families and of the special surfaces and we describe the Picard lattice of six surfaces with Picard-number 20 in the other families. Here by using lattice-theory and results on quadratic forms we compute the transcendental lattices of these surfaces. The methods are similar as the methods used by Barth in \cite{barth} for describing the K3-surfaces of \cite{basa}.

By a result of Morrison, \cite[Cor. 6.4]{Mo}, K3-surfaces with $\rho=19$ and $20$ have a Shioda-Inose structure, in particular this means that there is a birational map from the K3-surface to a Kummer surface. It is well known (cf. \cite{si}) that if $\rho=20$, then the Kummer surface is the quotient by an involution of a product-type abelian variety. When $\rho=19$ this is not always the case. In fact we use the transcendental lattices to show that in our cases the abelian variety is not a product of two elliptic curves. In this case we call the Shioda-Inose structure {\it simple}.\\ 
The paper is organized as follows: in section \ref{latticetheo} we recall some basic facts about lattices and quadratic forms and the construction of the families of K3-surfaces. Then section  \ref{transc} is entirely devoted to the computations of the transcendental lattices of the K3-surfaces of \cite{io}. In section \ref{abelian} we show that the Shioda-Inose structure of the surfaces with $\rho=19$ is simple. In section \ref{hessians} we compare our singular K3-surfaces with already  known surfaces, more precisely with the Hessians surfaces which are described in \cite{dard}: we see that all our singular K3-surfaces are Hessians of some cubic surface and we see that some of them are extremal elliptic K3-surfaces in the meaning of \cite{shizhang}. Finally in section \ref{conf} we recall the rational curves generating the Neron-Severi group of the K3-surfaces over $\mathbb{Q}$.\\
{\it I would like to thank Wolf Barth for letting me know about his paper \cite{barth} and for many discussions and  Slawomir Rams and Bert van Geemen for many useful comments.}
\section{Notations and preliminaries}\label{latticetheo}
\subsection{Lattices and quadratic forms}
A {\it lattice} L is a free $\Z$-module of finite rank with a $\Z$-valued symmetric bilinear form:
\begin{eqnarray*}
b:L\times L\longrightarrow \Z.
\end{eqnarray*}
An isomorphism of lattices preserving the bilinear form is called an {\it isometry}, $L$ is said to be {\it even} if the associate quadratic form to $b$ takes only even values, otherwise it is called {\it odd}. The {\it discriminant} $d(L)$ of $L$ is the determinant of the matrix of $b$, $L$ is said to be {\it unimodular} if $d(L)=\pm 1$. If $L$ is non-degenerate, i.e. $d(L)\not= 0$, then the {\it signature} of $L$ is a pair $(s_+, s_-)$ where $s_{\pm}$ denotes the multiplicity of the eigenvalue $\pm 1$ for the quadratic form on $L\otimes\R$, $L$ is called {\it positive-definite (negative-definite)} if the quadratic form associate to $b$ takes just positive (negative) values. We will denote by $U$ the {\it hyperbolic plane} i.e. a free $\Z$-module of rank 2 with bilinear form with matrix:
\begin{eqnarray*}
\left(\begin{array}{cc}
0&1\\
1&0\\
\end{array}\right)
\end{eqnarray*}
Moreover we denote by $E_8$ the unique even unimodular positive definite lattice of rank 8, with bilinear form with matrix:
\begin{eqnarray*}
\left(\begin{array}{cccccccc}
2&0&-1&0&0&0&0&0\\
0&2&0&-1&0&0&0&0\\
-1&0&2&-1&0&0&0&0\\
0&-1&-1&2&-1&0&0&0\\
0&0&0&-1&2&-1&0&0\\
0&0&0&0&-1&2&-1&0\\
0&0&0&0&0&-1&2&-1\\
0&0&0&0&0&0&-1&2\\
\end{array}\right)
\end{eqnarray*}
Let $L^{\vee}=Hom_{\Z}(L,\Z)=\{v\in L\otimes_{\mathbb{Z}}\mathbb{Q}~|~b(v,x)\in\mathbb{Z}~\mbox{for~all}~x\in L\}$ denotes the dual of the lattice $L$, then there is a natural embedding of $L$ in $L^{\vee}$ via $c\mapsto b(c,-)$, and we have:
\begin{lemma}\label{discr}(cf. \cite[Lemma 2.1, p. 12]{bpv})
If $L$ is a non-degenerate lattice with bilinear form $b$. Then\\
1. $[L^{\vee}:L]=|d(L)|$.\\
2. If $M$ is a submodule of $L$ with rank $M$=rank $L$, then
\begin{eqnarray*}
[L:M]^2=d(M)d(L)^{-1}.
\end{eqnarray*}
\end{lemma}
Let $A$ be a finite abelian group. A quadratic form on $A$ is a map:
\begin{eqnarray*}
q:A\longrightarrow \Q/2\Z
\end{eqnarray*}
together with a symmetric bilinear form:
\begin{eqnarray*}
b:A\times A \longrightarrow \Q/\Z
\end{eqnarray*}
such that:\\
1. $q(na)=n^2 q(a)$ for all $n\in\Z$ and $a\in A$\\
2. $q(a+a')-q(a)-q(a')\equiv 2 b(a,a')~(mod~2\Z)$\\
Let $L$ be a non-degenerate even lattice then the $\Q$-valued quadratic form on $L^{\vee}$ induces a quadratic form
\begin{eqnarray*}
q_{L}:L^{\vee}/L\longrightarrow \Q/2\Z
\end{eqnarray*}
called {\it discriminant-form} of $L$. By a result of Nikulin \cite[Cor. 1.9.4]{nikulin}, the signature and the discriminant form of an even lattice determines its {\it genus} (we do not need the exact definition here, cf. e.g. \cite{sloane}). \\
An embedding of lattices $M \hookrightarrow L$ is {\it primitive} if $L/M$ is free.
\begin{lemma}\label{discrform}(cf. \cite[Prop. 1.6.1]{nikulin})
Let $M \hookrightarrow L$ be a primitive embedding of non-degenerate even lattices and suppose $L$ unimodular then:\\
1. There is an isomorphism $M^{\vee}/M\cong (M^{\perp})^{\vee}/M^{\perp}$.\\
2. $q_{M^{\perp}}=-q_{M}$.
\end{lemma}

\bigskip

Let now X be an algebraic $K3$-surface, the group $H^2(X,\Z)$ with the intersection pairing has the structure of a lattice and by Poincar\'e duality it is unimodular. This is isometric to the {\it K3-lattice}:
\begin{eqnarray*}
\Lambda:=-E_8\oplus-E_8\oplus U\oplus U\oplus U
\end{eqnarray*}
(cf. \cite[Prop.3.2, p. 241]{bpv}).  The {\it Neron-Severi group} $NS(X)=H^2(X,\Z)\cap H^{1,1}(X)$ and its orthogonal complement $T_X$ in $H^2(X,\Z)$ (the {\it transcendental lattice}) are primitive sublattice of $H^2(X,\Z)$ and have signature $(1,\rho-1)$ and $(2,20-\rho)$, $\rho=$rank$(NS(X))$. By the Lemma \ref{discrform} we have\\
\begin{eqnarray*}
NS(X)^{\vee}/NS(X)\cong (T_X)^{\vee}/T_X
\end{eqnarray*}
and the discriminat-forms differ just by their sign. Moreover by the Lemma \ref{discr} we have $|NS(X)^{\vee}/NS(X)|=|(T_X)^{\vee}/T_X|=d(NS(X))$.

\bigskip

We recall some more facts  about $K3$-surfaces $X$ with $\rho=20$ ({\it singular} K3-surfaces, cf. \cite[p. 128]{si}). Denote by $\mathcal{Q}$  the set of
$2\times 2$ positive-definite even integral matrices:
\begin{eqnarray}\label{belle}
Q:=\left(\begin{array}{cc}
2a&c\\
c&2b\\
\end{array}\right),~~~~a,b,c\in\Z
\end{eqnarray}
with $d:=4ab-c^2> 0$ and $a,b >0$. 
We define $Q_1\sim Q_2$ if and only if $Q_1=~ ^t\gamma Q_2 \gamma$ for some $\gamma\in SL_2(\Z)$. Let $[Q]$ be the equivalence class of $Q$ and by $\mathcal{Q}/SL_2(\Z)$ the set of these equivalence classes. Then:
\begin{theorem}(cf. \cite[Thm. 4]{si}). 
The map $X\mapsto [T_X]$ estabilishes a bijective correspondence from the set of singular K3-surfaces onto $\mathcal{Q}/SL_2(\Z)$.
\end{theorem}
In particular $K3$-surfaces with $\rho=20$ are classified in terms of their transcendental lattice. By \cite[Thm. 2.3, p. 14]{buell}, we can assume that $Q$ is {\it reduced}, i.e. $-a\leq c\leq a\leq b$, and so $c^2\leq ab\leq d/3$. Recall the following:
\begin{theorem}(\cite[Theorem 2.4, p. 15]{buell})\label{bu}
With the exception of
\begin{eqnarray*}
1. \left(\begin{array}{cc}  2a&a\\a&2b\end{array}\right)\sim \left(\begin{array}{cc} 2a&-a\\-a&2b\end{array}\right); ~~~~~~~~~~~~~~~~~~~~~2. \left(\begin{array}{cc} 2a&b\\b&2a\end{array}\right)\sim \left(\begin{array}{cc} 2a&-b\\-b&2a \end{array}\right)
\end{eqnarray*}
no distinct reduced quadratic forms are equivalent.
\end{theorem}
Here the relation ``$\sim$'' is conjugation with a matrix of $SL_2(\mathbb{Z})$.\\
It is well known that the number of equivalence classes of forms of a given discriminat $d$, i.e. the {\it class number} of $d$, is finite. If there is only one class we say that $d$ has class number one. In some other cases we have one class per genus. In \cite[pp.81--82]{buell} with the assumption $g.c.d(a,c,b)=1$ all the discriminants of class number one and of one class per genus are listed. If $g.c.d(a,c,b)\not= 1$ then the form is a multiple of a primitive form.
\subsection{Families of K3-surfaces}
Let $G\subset SO(3)$ denotes the polyhedral group $T$, $O$ or $I$, and let  $\widetilde{G}\subset SU(2)$ be the corresponding binary groups. Let 
\begin{eqnarray*}
\sigma:SU(2)\times SU(2)\rightarrow SO(4, \R)
\end{eqnarray*}
denotes the classical $2:1$ covering. The images $\sigma(\widetilde{T}\times\widetilde{T}):=G_6$, $\sigma(\widetilde{O}\times\widetilde{O}):=G_8$  and  $\sigma(\widetilde{I}\times\widetilde{I}):=G_{12}$ in $SO(4,\R)$ are studied in \cite{io1}, where we show that there are  1-dimensional families in $\mathbb{P}_3(\mathbb{C})$ of $G_n$-invariant surfaces of degree $n$, which we denote by $X^n_{\lambda}$, $\lambda$ a parameter in $\mathbb{P}_1$. 
In \cite{basa} it is shown that the quotients $Y_{\lambda,G_n}$, $n=6,8,12$ are families of $K3$-surfaces where the general surface has Picard-number 19 and there are four surfaces with Picard-number 20. Then in \cite{io} by taking special normal subgroups of $G_n$ ($n=6,8$) and making the quotient of $X^6_{\lambda}$ resp. $X^8_{\lambda}$ by these subgroups we find six more pencils of K3-surfaces, using the notations there the subgroups are
\begin{eqnarray*}
\begin{array}{lllllll}
\mathcal{G}:~~~~T\times V&(TT)'&V\times V&O\times T&(OO)''&T\times T\\
\end{array}
\end{eqnarray*}
and the families of K3-surfaces are denoted by $Y_{\lambda,\mathcal{G}}$.
Here $V$ denotes the Klein four group in $SO(3, \R)$ and the groups $(TT)',(OO)''$ are described in \cite{io}, the others are the images in $SO(4,\R)$ of the direct product of binary subgroups of $SU(2)$. Moreover $T\times V$, $(TT)'$  are subgroups of index $3$ of $G_6$  and $V\times V$ has index $3$ in $T\times V$, $(TT)'$; $O\times T$, $(OO)''$ are subgroups of index $2$ of $G_8$ and $T\times T$ has index 2 in $O\times T$, $(OO)''$. 
In the families $Y_{\lambda,T\times V}$ and $Y_{\lambda,O\times T}$   the general surface has Picard-number 19 and we could identify four surfaces with Picard-number 20. We denote them by $Y_{\lambda,\mathcal{G}}^{(n,j)}$, where $n=6$, $\mathcal{G}=T\times V$ and $j=1,2,3,4$ or $n=8$, $\mathcal{G}=O\times T$ and $j=1,2,3,4$. In the other families we identify the Picard lattice of the following surfaces with $\rho=20$:
\begin{eqnarray*}
\begin{array}{l}
Y_{\lambda, (TT)'}^{(6,1)},~T_{\lambda, (TT)'}^{(6,2)},~Y_{\lambda, (OO)''}^{(8,1)},~Y_{\lambda, (OO)''}^{(8,4)}, Y_{\lambda, T\times T}^{(8,1)}, Y_{\lambda, T\times T}^{(8,4)}.
\end{array}
\end{eqnarray*}
We denote by $NS$ the Picard-lattice, by $T$ the transcendental lattice. We denote by $\Z_m(\alpha)$  the cyclic group $\Z_m$ with the quadratic form taking the value $\alpha\in\Q/2\Z$ on the generator of the group.
\section{Transcendental Lattices}\label{transc}
In this section we identify first the transcendental lattice of the singular $K3$-surfaces then of the surfaces with $\rho=19$. In each case we proceed as follows:\\
1. We determine generators for $NS^{\vee}/NS$ with the help of the intersection pairing $(-,-)$, which is defined on $NS$ (recall that $NS^{\vee}=\{v\in NS\otimes_{\mathbb{Z}}\mathbb{Q}~|~(v,x)\in\mathbb{Z}~\mbox{for~all}~x\in NS\}$).\\
2. We determine the discriminant-form of $NS$.\\
3. We use Lemma \ref{discrform} to determine the discriminant-form of $T$.\\
4. We list all the reduced quadratic forms which have the discriminant $d(T)=d(NS)$ (we will see that in each case the matrices have form 1 or 2 as in the Theorem \ref{bu}).\\
5. We use the discriminant form to determine $T$, in fact we see that when the rank is two the discriminants have class number one or one class per genus. When the rank is three in our cases the discriminants are {\it small}, Def. \ref{small}, and these have one class per genus. \\

\subsection{The singular cases}
{\it The family $Y_{\lambda, T\times V}$.}
We recall the following 3-divisible class of $NS$
\begin{eqnarray*}
\bar{L}'=L_1-L_2+L_4-L_5+L_1'-L_2'+L_4'-L_5'+L_1''-L_2''+L_4''-L_5''
\end{eqnarray*}
and the following 2-divisible classes of $NS$
\begin{eqnarray*}
h_1=L_1+L_3+L_5+L_1'+L_3'+L_5'+M_1+M_2,\\
h_2=L_1+L_3+L_5+L_1''+L_3''+L_5''+M_1+M_3.
\end{eqnarray*}The general $K3$-surface in the family has $\rho=19$ and the family contains four singular $K3$-surfaces.
The discriminant of the general $K3$-surface in the pencil is $2\cdot3\cdot 5$ which is the order of $NS^{\vee}/NS$ by the Lemma \ref{discr}. We specify the following generators:
\begin{eqnarray*}\label{5div}
\begin{array}{l}
M:=M_1+M_2+M_3/2,\\
N:=L_1-L_2+L_4-L_5-L_1'+L_2'-L_4'+L_5'/3,\\
L:=(3L_0-L_1-L_1'-L_1''-2L_2-2L_2'-2L_2''-3L_3-3L_3'-3L_3''\\
-2L_4-2L_4'-2L_4''-L_5-L_5'-L_5'')/5
\end{array}
\end{eqnarray*}
where
\begin{eqnarray*}
M^2=-3/2=1/2\mod 2\Z,\\
N^2=-8/3=4/3\mod 2\Z,\\
L^2=-18/5=2/5\mod 2\Z.
\end{eqnarray*}
Hence the dicriminant form of the Picard lattice is
\begin{eqnarray*}
\Z_2(1/2)\oplus\Z_3(4/3)\oplus\Z_5(2/5)\cong \Z_{30}(7/30)
\end{eqnarray*}
{\it The singular case $6,1(6,4)$}. Here the discriminant is $-3\cdot 5=-15$ and the generators of $NS^{\vee}/NS$ are $N$ and $L$. The dicriminant form is
\begin{eqnarray*}
\Z_3(4/3)\oplus\Z_5(2/5)=\Z_{15}(26/15)
\end{eqnarray*}
{\it The singular case $6,2(6,3)$}. Here the discriminant is $-2^2\cdot 3\cdot 5=-60$, and the generators are $M,N,L$ and another class $M'=M_4/2$ with $M'^2=-1/2=3/2$ $\mod 2\Z$. The discriminant form is
\begin{eqnarray*}
\Z_2(1/2)\oplus\Z_2(3/2)\oplus\Z_3(4/3)\oplus\Z_5(2/5)\cong\Z_2(1/2)\oplus\Z_{30}(97/30).
\end{eqnarray*}
The discriminant form of the transcendental lattice differs by the previous form just by the sign, hence in the general case is
\begin{eqnarray*}
\Z_{30}(53/30)
\end{eqnarray*}
and in the special cases is
\begin{eqnarray*}
\begin{array}{ll}
6,1~(6,4):&\Z_{15}(4/15),\\
6,2~(6,3):&\Z_{2}(3/2)\oplus\Z_{30}(23/30).\\
\end{array}
\end{eqnarray*}
Here we identify the transcendental lattices of these four singular $K3$-surfaces, and in the next section of the general $K3$-surface.\\
{\it The singular case $6,1~(6,4)$}. We classify all the reduced matrices with discriminant $15$ (one representant per class, cf. \cite[pp.19--20]{buell}).
We have just the following possibilities
\begin{eqnarray*}
\left(\begin{array}{cc}
2&1\\
1&8\\
\end{array}\right), ~~~~A:=\left(\begin{array}{cc}
4&1\\
1&4\\
\end{array}\right).
\end{eqnarray*}
By taking the generator $(4/15, -1/15)$ and the bilinear form defined by $A$, we find a lattice $\Z_{15}(4/15)$ which is exactly the lattice $T^{\vee}/T$ hence $T=A$.\\
{\it The singular case $6,2~(6,3)$}. We classify all the reduced matrices with discriminant $60$( cf. \cite[pp.19--20]{buell}).
We have just the following possibilities
\begin{eqnarray*}
\left(\begin{array}{cc}
2&0\\
0&30\\
\end{array}\right), ~~~~B:=\left(\begin{array}{cc}
6&0\\
0&10\\
\end{array}\right), ~~~~
\left(\begin{array}{cc}
4&2\\
2&16\\
\end{array}\right), ~~~~\left(\begin{array}{cc}
8&2\\
2&8\\
\end{array}\right).
\end{eqnarray*}
By taking the generators $(1/2, 0)$ and  $(1/3, 1/10)$ and the quadratic form $B$ we find a lattice $\Z_{2}(3/2)\oplus\Z_{30}(23/30)$ which is exactly the lattice $T^{\vee}/T$, hence $T=B$.\\





{\it The family $Y_{\lambda,(TT)'}$}. We recall the following 3-divisible class in $NS$:
\begin{eqnarray*}
\bar{L}=N_1-N_2+N_3-N_4+N_5-N_6+N_7-N_8+N_9-N_{10}+N_{11}-N_{12}.
\end{eqnarray*}
Now we identify the transcendental lattice of $Y_{\lambda,(TT)'}^{(6,1)}$ and of  $Y_{\lambda,(TT)'}^{(6,2)}$.\\
{\it The singular case $6,1$.} In this case the discriminant is $-3\cdot 5=-15$ and we have the following generators of $NS^{\vee}/NS$:
\begin{eqnarray*}
\begin{array}{l}
N:=(N_1-N_2+N_3-N_4-N_5+N_6-N_7+N_8)/3,\\
L:=(3L_3-3L_3')/5,\\
\end{array}
\end{eqnarray*}
where
\begin{eqnarray*}
N^2=-8/3=4/3 \mod 2\Z,\\
L^2=-18/5=2/5 \mod 2\Z.
\end{eqnarray*}
Hence the transcendental lattice is the same as in the case of $Y_{\lambda, T\times V}^{(6,1)}$.\\
{\it The singular case $6,2$.} Recall the following 2-divisible classes in $NS$:
\begin{eqnarray*}
N_1+C_1+N_4+N_5+C_2+N_8+M_1+M_2,\\
N_1+C_1+N_4+N_9+C_3+N_{12}+M_1+M_3.\\
\end{eqnarray*}
The discriminant is $-2^2\cdot 3\cdot 5=-60$ and the classes
\begin{eqnarray*}
N,~~M=M_1+M_2+M_3/2,~~M'=N_5+C_2+N_8+M_1+M_3/2,~~L
\end{eqnarray*}
are generators for $NS^{\vee}/NS$. Where
\begin{eqnarray*}
\begin{array}{l}
N^2=-8/3=4/3\mod 2\Z,\\
M^2=1/2\mod 2\Z,\\
M'^2=3/2\mod 2\Z,\\
L^2=-18/5=2/5\mod 2\Z.\\
\end{array}
\end{eqnarray*}
Hence the transcendental lattice is the same as in the case of $Y_{\lambda, T\times V}^{(6,2)}$.\\
{\it The family $Y_{\lambda,O\times T}$.} Recall the following 2-divisible class of $NS$:
\begin{eqnarray*}
\bar{L}'=L_1+L_3+L_5+L_1'+L_3'+L_5'+M_1+M_2,
\end{eqnarray*}
and the following 3-divisible class of $NS$:
\begin{eqnarray*}
k_1=L_1-L_2+L_4-L_5-L_1'+L_2'-L_4'+L_5'+N_1-N_2+N_3-N_4.
\end{eqnarray*}
The general surface in the pencil has $\rho=19$ and we have four surfaces with $\rho=20$.
The discriminant of the general $K3$-surface in the pencil is $2^3\cdot 3\cdot 7=168$. We specify the following generators of $NS^{\vee}/NS$:
\begin{eqnarray*}\label{7div}
\begin{array}{l}
M:=L_1+L_3+L_5+M_2/2,\\
M':=L_1+L_3+L_5+M_3/2,\\
R:=R_2/2,\\
N:=N_1-N_2-N_3+N_4/3,\\
L:=(2L_2''+4L_0-2L_1-2L_1'+3L_2+3L_2'-3L_3-3L_3'-2L_4-2L_4'-L_5-L_5')/7\\
\end{array}
\end{eqnarray*}
where
\begin{eqnarray*}
\begin{array}{l}
M^2=-2=0 \mod 2\Z,\\
M'^2=-2=0\mod 2\Z,\\
R^2=-1/2=3/2\mod 2\Z,\\
N^2=-4/3=2/3\mod 2\Z,\\
L^2=-16/7=12/7\mod 2\Z.\\
\end{array}
\end{eqnarray*}
Observe that the classes $M$, $M'$ and $L$ are not orthogonal to eachother in fact $M\cdot M'=1/2~\mod 2\Z$ and $M\cdot L=M'\cdot L=1~\mod 2\Z$.
Hence the discriminant form of the Picard lattice is:
\begin{eqnarray*}
\Z_2(0)\oplus\Z_2(0)\oplus\Z_2(3/2)\oplus\Z_3(2/3)\oplus\Z_7(12/7)) \cong \Z_2(0)\oplus\Z_2(0)\oplus\Z_{42}(79/42).
\end{eqnarray*}
{\it The singular case $8,1$.} Here the discriminant is $-2^2\cdot 7=-28$ and the generators for $NS^{\vee}/NS$ are $M$, $M'$ and $L$. The discriminant form is
\begin{eqnarray*}
\Z_2(0)\oplus\Z_2(0)\oplus\Z_7(12/7)) \cong\Z_2(0)\oplus\Z_{14}(12/7).
\end{eqnarray*}
{\it The singular case $8,2$.} The discriminant is $-2^2\cdot 3\cdot 7=-84$ and the generators for $NS^{\vee}/NS$ are $M+R$, $M'+R$, $N$  and $L$.  The discriminant form is
\begin{eqnarray*}
\Z_2(3/2)\oplus\Z_2(3/2)\oplus\Z_3(2/3))\oplus\Z_7(12/7) \cong\Z_2(3/2)\oplus\Z_{42}(163/42)=\Z_2(3/2)\oplus\Z_{42}(79/42).
\end{eqnarray*}
{\it The singular case $8,3$.} Here the discriminant is $-2^3\cdot 3\cdot 7=-168$ and the generators for $NS^{\vee}/NS$ are $R$,
\begin{eqnarray*}
R'=M_1+2C+3M_2/4,
\end{eqnarray*}
$N$  and $L$, where $R'^2=1/4\mod 2\Z$. The discriminant form is
\begin{eqnarray*}
\Z_2(3/2)\oplus\Z_4(1/4)\oplus\Z_3(2/3)\oplus\Z_7(12/7)) \cong\Z_2(3/2)\oplus\Z_{84}(221/84)=\Z_2(3/2)\oplus\Z_{84}(53/84).
\end{eqnarray*}
{\it The singular case $8,4$.} Recall the 2-divisible class in $NS$
\begin{eqnarray*}
L_1+L_3+L_5+N_1+C+N_4+R_2+M_1
\end{eqnarray*}
The discriminant is $-2^2\cdot 7=-28$ and the generators  for $NS^{\vee}/NS$ are $L'+R$,
\begin{eqnarray*}
M''=M_1+M_2+R_2/2,
\end{eqnarray*}
and $L$, where $M''^2=1/2 \mod 2\Z$.\\
The discriminant form is
\begin{eqnarray*}
\Z_2(3/2)\oplus\Z_2(1/2)\oplus\Z_7(12/7)) \cong\Z_2(3/2)\oplus\Z_{14}(31/14)\cong\Z_2(3/2)\oplus\Z_{14}(3/14) (mod~2\Z).
\end{eqnarray*}
The discriminant of the transcendental lattice differs by the previous form just by the sign, hence in the general case is
\begin{eqnarray*}
\Z_2(0)\oplus\Z_2(0)\oplus\Z_{42}(5/42)
\end{eqnarray*}
and in the special cases is
\begin{eqnarray*}
\begin{array}{ll}
8,1:&\Z_2(0)\oplus\Z_{14}(2/7),\\
8,2:&\Z_2(1/2)\oplus\Z_{42}(5/42),\\
8,3:&\Z_2(1/2)\oplus\Z_{84}(115/84),\\
8,4:&\Z_2(1/2)\oplus\Z_{14}(25/14).\\
\end{array}
\end{eqnarray*}
Here we identify the transcendental lattice for this four singular cases, and in the next section for the general $K3$-surface.\\
{\it The singular case $8,1$.} We classify all the reduced matrices with discriminant $28$ (\cite[pp.19--20]{buell}). We have just the following possibilities:
\begin{eqnarray*}\label{matriceB}
A:=\left(\begin{array}{cc}
2&0\\
0&14\\
\end{array}\right), ~~~~B:=\left(\begin{array}{cc}
4&2\\
2&8\\
\end{array}\right).
\end{eqnarray*}
Now take the form $B$ and the generators $(0,1/2)$ and $(3/14,1/14)$. These span exactly the lattice we were looking for.\\
{\it The singular case $8,2$.}We classify all the reduced matrices with discriminant $84$ (\cite[pp.19--20]{buell}). We have the following four cases:
\begin{eqnarray*}
\left(\begin{array}{cc}
2&0\\
0&42\\
\end{array}\right), ~~~~\left(\begin{array}{cc}
6&0\\
0&14\\
\end{array}\right),~~~~
\left(\begin{array}{cc}
4&2\\
2&22\\
\end{array}\right), ~~~~C:=\left(\begin{array}{cc}
10&4\\
4&10\\
\end{array}\right).
\end{eqnarray*}
Now we take the form $C$ and the generators $(1/2, 0)$ and $(8/21, -19/42)$ and we are done.\\
{\it The singular case $8,3$.}We classify all the reduced matrices with discriminant $168$ (\cite[pp.19--20]{buell}). We have the following four cases
\begin{eqnarray*}
\left(\begin{array}{cc}
2&0\\
0&84\\
\end{array}\right), ~~~~\left(\begin{array}{cc}
6&0\\
0&28\\
\end{array}\right),~~~~
E:=\left(\begin{array}{cc}
12&0\\
0&14\\
\end{array}\right), ~~~~\left(\begin{array}{cc}
4&0\\
0&42\\
\end{array}\right).
\end{eqnarray*}
Now we take the form $E$ and the generators $(1/2, 1/2)$ and $(1/12,1/7)$. These span exactly the lattice we were looking for.\\
{\it The singular case $8,4$.}
The discriminant is $28$ like in the case of $8,1$. Now by taking the form $A$ and the generators $(1/2,0)$ and $(0,5/14)$ we are done.\\
{\it The family $Y_{\lambda,(OO)''}$.} Recall the following 2-divisible class of $NS$
\begin{eqnarray*}
\bar{L}=M_1+M_2+M_3+M_4+R_1+R_3+R_1'+R_3'
\end{eqnarray*}
We identify the transcendental lattices of the surfaces $Y_{\lambda,(OO)''}^{(8,1)}$ and $Y_{\lambda,(OO)''}^{(8,4)}$ .\\
{\it The singular case $8,1$}. 
In this case the the discriminant is $-2^2\cdot 7=-28$ and we have the following generators in $NS^{\vee}/NS$
\begin{eqnarray*}
\begin{array}{l}
L:=2L_2+4L_4-2L_2'-4L_4'/7,\\
M:=R_1+R_3+M_1+M_3/2,\\
M':=R_1+R_3+M_1+M_4/2,\\
\end{array}
\end{eqnarray*}
where
\begin{eqnarray*}
\begin{array}{l}
L^2=12/7 \mod 2\Z,\\
M^2=M'^2=0 \mod 2\Z.\\
\end{array}
\end{eqnarray*}
Hence the transcendental lattice is the same as in the case $Y_{\lambda,O\times T}^{(8,1)}$.\\
{\it The singular case $8,4$.}
Recall the following 4-divisible class in NS\\
\begin{eqnarray*}
W:=R_1+2R_2+3R_3+R_1'+2R_2'+3R_3'+2N_1+2C_1+3M_1+M_2+2N_3+2C_2+3M_3+M_4.
\end{eqnarray*}
Moreover specify the classes:
\begin{eqnarray*}
\begin{array}{l}
v_1:=R_1+2R_2+3R_3/4,\\
v_2:=R_1'+2R_2'+3R_3'/4,\\
v_3:=2N_1+2C_1+3M_1+M_2/4,\\
v_4:=2N_3+2C_2+3M_3+M_4/4,.\\
\end{array}
\end{eqnarray*}
The discriminant is $-2^4\cdot 7=-112$ and the generators of $NS^{\vee}/NS$ are
\begin{eqnarray*}
v_1+v_3/4,~~~~~v_2+v_4/4,~~~~~L
\end{eqnarray*}
with
\begin{eqnarray*}
(v_1+v_3/4)^2=(v_2+v_4/4)^2=0~~~\mod 2\Z.
\end{eqnarray*}
The discriminant form of the Picard lattice is
\begin{eqnarray*}
\Z_4(0)\oplus\Z_4(0)\oplus\Z_7(12/7)=\Z_4(0)\oplus\Z_{28}(12/7).
\end{eqnarray*}
Hence the discriminant form of the transcendental lattice is
\begin{eqnarray*}
\Z_4(0)\oplus\Z_{28}(2/7)
\end{eqnarray*}
We classify all the reduced matrices with discriminant $112$, these are
\begin{eqnarray*}
\left(\begin{array}{cc}
2&0\\
0&56\\
\end{array}\right), ~~~~\left(\begin{array}{cc}
4&0\\
0&28\\
\end{array}\right),~~~~
F:=\left(\begin{array}{cc}
8&0\\
0&14\\
\end{array}\right), ~~~~\left(\begin{array}{cc}
8&4\\
4&16\\
\end{array}\right).
\end{eqnarray*}
We take the matrix $F$ and the generators $(1/4, 1/2)$ and $(1/4, 9/14)$, so we are done.\\
{\it The family $Y_{\lambda, T\times T}$.} A similar computation as before shows that in the singular case $8,1$, resp. $8,4$  the transcendental lattice has bilinear form with intersection matrix: 
$$
\left(\begin{array}{cc}
2&1\\
1&4\\
\end{array}\right),~~~~~~~~~~~~~~~~\mbox{resp}~~~~~~~~~~~~~~~~~~~~~~~~~~~~\left(\begin{array}{cc}
4&2\\
2&8\\
\end{array}\right).
$$
\begin{rem}
Observe that if the reduced matrices had not been as in case 1 or 2 of Theorem \ref{bu} we would find two different isomorphism classes of K3-surfaces with the same discriminant and the same discriminant form (cf. \cite{shizhang} p. 3). 
\end{rem}
\subsection{The general cases} Here we identify the transcendental lattice of the general surfaces, $\rho=19$ in the families $Y_{\lambda, T\times V}$ and $Y_{\lambda,O\times T}$. In the last section we have identified the discriminant form of the transcendental lattice, we use it to determine $T$. We need the following: \\
\begin{defi}\label{small} (cf.\cite[Def. 1.1]{barth}) The discriminant $d=d_{NS}=-d_T$ is {\rm small} if $4\cdot d$ is not divisible by $k^3$ for any non square natural number $k$ congruent to $0$ or $1$ modulo 4.
\end{defi}
Then if $d_T$ is small , the lattice $T$ is uniquely determined by its genus (cf. \cite[Thm. 21, p. 395]{sloane}), hence by signature and discriminant form.\\
{\it The family $Y_{\lambda, T\times V}$}. The candidate lattice is
\begin{eqnarray*}
T_0:=\left(\begin{array}{rrr}
4&1&0\\
1&4&0\\
0&0&-2\\
\end{array}
\right)
\end{eqnarray*}
this has discriminant -30, and taking the generator
\begin{eqnarray*}
f_1:=\left(\begin{array}{r}
4/15\\
-1/15\\
1/2\\
\end{array}
\right)
\end{eqnarray*}
one computes $q_{T_0}(f_1)=-7/30=53/30 \mod 2\Z$, hence the discriminant form is $\Z_{30}(53/30)$.  Since $d_T=-30$ is small the transcendental lattice of the general K3-surface is $T_0$.\\
{\it The family $Y_{\lambda,O\times T}$}. The candidate lattice is
\begin{eqnarray*}
T_1:=\left(\begin{array}{rrr}
10&4&0\\
4&10&0\\
0&0&-2\\
\end{array}
\right)
\end{eqnarray*}
this has discriminant -168, and taking the generators
\begin{eqnarray*}
f_1:=\left(\begin{array}{r}
1/2\\
0\\
1/2\\
\end{array}
\right),~~~~~
f_2:=\left(\begin{array}{r}
1/2\\
0\\
-1/2\\
\end{array}
\right),~~~~~
f_3:=\left(\begin{array}{r}
8/21\\
-19/42\\
0\\
\end{array}
\right)
\end{eqnarray*}
we find $q_{T_1}(f_i)=0 \mod 2\Z$, $i=1,2$ and $q_{T_1}(f_3)=5/42 \mod 2\Z$, hence the discriminant form is $\Z_{2}(0)\oplus\Z_{2}(0)\oplus \Z_{42}(5/42)$. Since $d_T=-168$ is small we have $T=T_1$.\\
We collect the results in the table \ref{tab1}. We recall also the results of \cite{barth} about the general surfaces of the families $Y_6(\lambda)$,   $Y_8(\lambda)$, $Y_{12}(\lambda)$ and also about the singular surfaces in these pencils, Barth computed the transcendental lattices of the singular surfaces too, but he did not published his result. In the table we write also the discriminants of the lattices.

\begin{footnotesize}
\begin{table}\caption{Transcendental Lattices}\label{tab1}
\begin{tabular}{l|c|cccc}
Family&general surface&\multicolumn{4}{c}{singular surfaces}\\
\hline
$Y_{\lambda, G_6}$& $\left(\begin{array}{rrr}2&1&0\\1&8&0\\0&0&-6\end{array}\right)$& $\left(\begin{array}{rr}2&1\\1&8\end{array}\right)$&
$\left(\begin{array}{rr}2&0\\0&30\end{array}\right)$&
$\left(\begin{array}{rr}2&0\\0&30\end{array}\right)$&
$\left(\begin{array}{rr}2&1\\1&8\end{array}\right)$\\
\hline
$d$&$-90$&$15$&$60$&$60$&$15$\\
\hline
\hline
$Y_{\lambda, G_8}$&$\left(\begin{array}{rrr}6&0&0\\0&28&0\\0&0&-2\end{array}\right)$&
$\left(\begin{array}{rr}2&0\\0&14\end{array}\right)$&
$\left(\begin{array}{rr}6&0\\0&14\end{array}\right)$&
$\left(\begin{array}{rr}6&0\\0&28\end{array}\right)$&
$\left(\begin{array}{rr}4&0\\0&28\end{array}\right)$\\
\hline
$d$&$-336$&$28$&$84$&$168$&$112$\\
\hline
\hline
$Y_{\lambda, G_{12}}$&$\left(\begin{array}{rrr}4&2&0\\2&34&0\\0&0&-30\end{array}\right)$&
$\left(\begin{array}{rr}12&6\\6&58\end{array}\right)$&
$\left(\begin{array}{rr}6&0\\0&220\end{array}\right)$&
$\left(\begin{array}{rr}6&0\\0&132\end{array}\right)$&
$\left(\begin{array}{rr}4&2\\2&34\end{array}\right)$\\
\hline
$d$&$-3960$&$660$&$1320$&$792$&$132$\\
\hline
\hline
$Y_{\lambda, T\times V}$&$\left(\begin{array}{rrr}4&1&0\\1&4&0\\0&0&-2\end{array}\right)$&
$\left(\begin{array}{rr}4&1\\1&4\end{array}\right)$&
$\left(\begin{array}{rr}6&0\\0&10\end{array}\right)$&
$\left(\begin{array}{rr}6&0\\0&10\end{array}\right)$&
$\left(\begin{array}{rr}4&1\\1&4\end{array}\right)$\\
\hline
$d$&$-30$&$15$&$60$&$60$&$15$\\
\hline
\hline
$Y_{\lambda,O\times T}$&$\left(\begin{array}{rrr}10&4&0\\4&10&0\\0&0&-2\end{array}\right)$&
$\left(\begin{array}{rr}4&2\\2&8\end{array}\right)$&
$\left(\begin{array}{rr}10&4\\4&10\end{array}\right)$&
$\left(\begin{array}{rr}12&0\\0&14\end{array}\right)$&
$\left(\begin{array}{rr}2&0\\0&14\end{array}\right)$\\
\hline
$d$&$-168$&$28$&$84$&$168$&$28$\\
\hline
\hline
$Y_{\lambda,(TT)'}$&-&$\left(\begin{array}{rr}4&1\\1&4\end{array}\right)$&
$\left(\begin{array}{rr}6&0\\0&10\end{array}\right)$&-&-\\
\hline
$d$&-&$15$&$60$&-&-\\
\hline
\hline
$Y_{\lambda,(OO)''}$&-&$\left(\begin{array}{rr}2&0\\0&14\end{array}\right)$&-&-&$\left(\begin{array}{rr}8&0\\0&14\end{array}\right)$\\
\hline
$d$&-&$28$&-&-&$112$\\
\hline
\hline
$Y_{\lambda, T\times T}$&-&$\left(\begin{array}{rr}2&1\\1&4\end{array}\right)$&-&-&$\left(\begin{array}{rr}4&2\\2&8\end{array}\right)$\\
\hline
$d$&-&$7$&-&-&$28$
\end{tabular}
\end{table}
\end{footnotesize}

\subsection{Moduli curve}\label{iso}
 Let
\begin{eqnarray*}
\Omega=\{[\omega]\in\mathbb{P}(\Lambda\otimes\mathbb{C})|(\omega,\omega)=0; (\omega,\bar{\omega})>0\},
\end{eqnarray*}
this is an open subset in a quadric in $\mathbb{P}^{21}$. If $X$ is a K3-surface and $\omega_X\in H^{2,0}(X)$, then it is well known that $\omega_X\in\Omega$ and it is called a {\it period point}. Moreover also the converse is true: each point of $\Omega$ occurs as period point of some K3-surface, this is the so called {\it surjectivity of the period map} (cf. \cite[Thm. 14.2]{bpv}). Now let $M\subset\Lambda$ be a sublattice of signature $(1,\rho-1)$ and define:
\begin{eqnarray*}
 \Omega_M=\{[\omega]\in \Omega |(\omega,\mu)=0~\mbox{for~all}~\mu\in M\}.
\end{eqnarray*}
This has dimension $20-\rho=20-$rank $M$. If rank $M$=19 then this space is a curve. Let $X$ be a K3-surface with $\rho=19$ since in this case the embedding of $T_X$ in $\Lambda$ is unique up to isometry of $\Lambda$ (cf. \cite[Cor. 2.10]{Mo}), $T_X$ determines  $\Omega_M$, with $M=T_X^{\perp}=NS(X)$ and so the moduli curve, which classify the K3-surfaces. Hence in our cases the transcendental lattices given in the table \ref{tab1} identify the moduli curve of the K3-surfaces in the families $Y_{\lambda, T\times V}$ and $Y_{\lambda,O\times T}$ (in the case $\rho=19$).

\section{Shioda-Inose structure}\label{abelian}
By a result of Morrison $K3$-surfaces with $\rho=19$ or $\rho=20$ admit a Shioda-Inose structure. Before discussing our cases we recall some facts.
\begin{defi}(cf. \cite[Def. 6.1]{Mo})
A K3-surface $X$ admits a {\it Shioda-Inose structure} if there is a Nikulin Involution $\iota$ on $X$ with rational quotient map 
$\pi:X--\rightarrow Y$ such that $Y$ is a Kummer surface, and $\pi_*$ induces an Hodge isometry $T_X(2)\cong T_Y$.
\end{defi}
Hence we have the following diagram:
\begin{equation} 
\xymatrix{
      & A \ar@{->}[ld]\ar@{-->}[rd]&                   & 
X   \ar@{->}[rd] \ar@{-->}[ld]& & \\
A/i &\ar[l] &Y&\ar[r]&                    
                        X/\iota      
} 
\end{equation}
where $A$ is the complex torus whose Kummer-surface is $Y$, $\iota$ is a Nikulin involution, i.e. an involution with $8$ fix-points on $X$, $i$ is an involution on $A$ with 16 fix-points  
and the rational maps from $A$ to $Y$ and from $X$ to $Y$ are 2:1. By definition we have $T_X(2)\cong T_Y$ and by \cite[Prop. 4.3]{Mo}, we have $T_A(2)\cong T_Y$ hence the diagram induces an Hodge isometry $T_X\cong T_A$.\\
In our cases the $K3$-surface  which we consider are algebraic hence $A$ is an abelian variety (cf. \cite [Thm. 6.3, (ii)]{Mo}). Moreover whenever $X$ is an algebraic $K3$-surface and $\rho(X)=19$ or $20$ then $X$ admits always a Shioda-Inose structure (cf. \cite[Cor. 6.4]{Mo}). Whenever $\rho=20$ Shioda and Inose show that $A=C_1\times C_2$ where $C_1$ and $C_2$ are elliptic curves
\begin{eqnarray*}
C_i=\ci/\Z+\Z\cdot \tau_i,~~~~~~~~~~~~~~~~~~~~~~~~~~~~~~~~~~i=1,2
\end{eqnarray*}
whith 
\begin{eqnarray*}
\tau_1=(-c+\sqrt{-d})/2a,~~~~~~~~~~~\tau_2=(c+\sqrt{-d})/2,~~~~~~~~~~~~~~~~(d=4ab-c^2)
\end{eqnarray*}

We show that in the case of the general $K3$-surfaces of the families $Y_{\lambda, T\times V}$ and of $Y_{\lambda,O\times T}$ the abelian surface $A(\lambda)$ is {\it simple}, i.e. it is not a product of elliptic curves, in this case we say that the Shioda-Inose structure is {\it simple}.\\
The transcendental lattice $T_{A(\lambda)}$ has rank 3 hence its orthogonal complement $NS_{A(\lambda)}$ in $U^3$ has rank 3 too and we have $NS(A(\lambda))\cong T(Y(\lambda))(-1)$ because by \cite[Thm. 21, p. 395]{sloane}, the lattices are uniquely determined. We use this fact to show:\\
\begin{theorem}
For general $\lambda$, $A=A(\lambda)$ is not a product of elliptic curves. 
\end{theorem}
\bprf (cf.\cite[ Thm. 5.1]{barth})
We show that $A$ does not contain any elliptic curve $C$, i.e. a curve with $C^2=0$.\\
{\it The general surface in $Y_{\lambda, T\times V}$}: We have intersection form on the transcendental lattice with matrix
\begin{eqnarray*}
T_0:=\left(\begin{array}{rrr}
4&1&0\\
1&4&0\\
0&0&-2\\
\end{array}
\right)
\end{eqnarray*}
hence the form on $NS_A$ is
\begin{eqnarray*}
\left(\begin{array}{rrr}
-4&-1&0\\
-1&-4&0\\
0&0&2\\
\end{array}
\right).
\end{eqnarray*}
The associated quadratic form is 
\begin{eqnarray*}
-4x^2-2xy-4y^2+2z^2,~~~~x,y,z\in\Z.
\end{eqnarray*}
If $A$ contains an elliptic curve, then there are $x,y,z\in\Z$ with
\begin{eqnarray*}
2z^2=4x^2+2xy+4y^2
\end{eqnarray*}
hence
\begin{eqnarray*}
8z^2=16x^2+8xy+16y^2
\end{eqnarray*}
Put $u=4x+y$, then
\begin{eqnarray}\label{laprima}
8z^2=u^2+15 y^2.
\end{eqnarray}
Hence we have $u^2=3z^2$  $\mod 5\Z $, since 3 is not a square modulo 5 we have $u=z=0$ $\mod 5\Z$, hence $u=5u_1$, $z=5z_1$, so
\begin{eqnarray}\label{laseconda}
3y^2=5(8z^2- u^2)
\end{eqnarray}
hence $y=5 y_1$ and substituting in (\ref{laseconda}) and dividing by 5 we find
\begin{eqnarray*}
15 y_1^2=8z^2-u^2
\end{eqnarray*}
which is the same as (\ref{laprima}).\\
{\it The general surface in $Y_{\lambda,O\times T}$}: We have intersection form on the transcendental lattice with matrix:
\begin{eqnarray*}
T_1:=\left(\begin{array}{rrr}
10&4&0\\
4&10&0\\
0&0&-2\\
\end{array}
\right).
\end{eqnarray*}
Hence the form on $NS_A$ is 
\begin{eqnarray*}
\left(\begin{array}{rrr}
-10&-4&0\\
-4&-10&0\\
0&0&2\\
\end{array}
\right).
\end{eqnarray*}
The quadratic form is
\begin{eqnarray*}
-10x^2-8xy-10y^2=2z^2,~~~~x,y,z\in\Z
\end{eqnarray*}
If $A$  contains an elliptic curve, then there are $x,y,z\in\Z$ with
\begin{eqnarray*}
2z^2=10x^2+8xy+10y^2
\end{eqnarray*}
hence dividing by 2 and multiplying by 5 we find
\begin{eqnarray*}
5z^2=25x^2+20xy+25y^2=(5x+2y)^2+21y^2.
\end{eqnarray*}
Put $u=5x+2y$, then
\begin{eqnarray}\label{laterza}
5z^2=u^2+21 y^2.
\end{eqnarray}
Hence we have $u^2=5z^2$ $\mod 7\Z$. Since 5 is not a square modulo 7 we have  $u=7u_1$, $z=7z_1$, so we obtain
\begin{eqnarray}\label{laquarta}
3y^2=7(5z_1^2- u_1^2)
\end{eqnarray}
hence $3y^2=0$ $\mod 7\Z$. Since 3 is not a square modulo 7 we have $y=7y_1$  and substituting in (\ref{laquarta}) and dividing by 7 we find
\begin{eqnarray*}
21 y_1^2=5z_1^2-u_1^2
\end{eqnarray*}
which is again (\ref{laterza}).\eprf
\section{Hessians and extremal elliptic K3-surfaces}\label{hessians}
Many of the singular K3-surfaces of this article appear already in other realizations.\\ 
In \cite{dard} Dardanelli and van Geemen give a criteria to estabilish if a singular K3-surface is the desingularization of the Hessian of a cubic surface:
\begin{prop}(cf. \cite[Prop. 2.4.1]{dard})
Let T be an even lattice of rank 2,
\begin{eqnarray*}
T=\left(\begin{array}{cc}
2n&a\\
a&2m
\end{array}\right).
\end{eqnarray*}
There is a primitive embedding $T\hookrightarrow T_{Hess}$ if and only if at least one among a, n and m is even. In this case T embeds in $U \oplus U(2)$.
\end{prop}
Here  $T_{Hess}=U\oplus U(2)\oplus A_2(-2)$. If we look in table \ref{tab1} we see that all our singular $K3$-surfaces  are desingularizations of Hessians of cubic surfaces. In particular Dardanelli and van Geemen study explicitely the singular K3-surfaces with
 \begin{eqnarray*}
T=\left(\begin{array}{cc}
4&1\\
1&4
\end{array}\right).
\end{eqnarray*}
 They call the surface $X_{10}$ and show that it is the desingularization of the Hessian of the cubic surface with 10 Eckardt points. The latter has e.g. the following equation in $\mathbb{P}^4$
\begin{eqnarray*}
\sum_{i=0}^{4}x_i^3=0,~~~~~~\sum_{i=0}^{4}x_i=0.
\end{eqnarray*}
Finally observe that the singular surfaces of the families $Y_{\lambda,G_6}$, $Y_{\lambda, T\times V}$ and $Y_{\lambda,(TT)'}$ are {\it extremal elliptic K3-surfaces}, in the sense of Shimada and Zhang (cf. \cite{shizhang}), in fact these are the numbers: 322, 173, 102, 148, 276 in their list in  \cite[Table 2, pp. 15-24]{shizhang}.\\

\section{Figures: Configurations of rational curves}\label{conf}
In this section we recall the configurations of $(-2)$-rational curves generating the Neron-Severi group over $\mathbb{Q}$. In the case of the families  $Y_{\lambda, T\times V}$ and $Y_{\lambda, O\times T}$ the curves $L_i$, $L_i'$ and $L_i''$ on the general K3-surface are also contained in the Neron-Severi group of the singular K3-surfaces, but we do not draw their configuration again. Moreover since the singular surfaces $Y_{\lambda, T\times V}^{(6,1)}$ and $Y_{\lambda, T\times V}^{(6,4)}$, as the surfaces $Y_{\lambda, T\times V}^{(6,2)}$ and $Y_{\lambda, T\times V}^{(6,3)}$ have the same graph, we draw just one picture.

\bigskip




   \begin{center}
   \begin{psfrags}
     \psfrag{YTT}{$Y_{\lambda,G_6}$}
     \psfrag{YOO}{$Y_{\lambda,G_8}$}
     \psfrag{L1}{$L_1$}
     \psfrag{L2}{$L_2$}
     \psfrag{L3}{$L_3$}
     \psfrag{L4}{$L_4$}
     \psfrag{L5}{$L_5$}
\psfrag{L1'}{$L_1'$}
     \psfrag{L2'}{$L_2'$}
     \psfrag{L3'}{$L_3'$}
     \psfrag{L4'}{$L_4'$}
     \psfrag{L5'}{$L_5'$}
     \psfrag{L2''}{$L_2''$}
\psfrag{M1}{$M_1$}
\psfrag{M2}{$M_2$}
\psfrag{M3}{$M_3$}
\psfrag{M4}{$M_4$}
 \psfrag{L0}{$L_0$}
\psfrag{N1}{$N_1$}
\psfrag{N2}{$N_2$}
\psfrag{N3}{$N_3$}
\psfrag{N4}{$N_4$}
\psfrag{N5}{$N_5$}
\psfrag{N6}{$N_6$}
\psfrag{N7}{$N_7$}
\psfrag{N8}{$N_8$}
\psfrag{R1}{$R_1$}
\psfrag{R2}{$R_2$}
\psfrag{R3}{$R_3$}
\psfrag{R1'}{$R_1'$}
\psfrag{R2'}{$R_2'$}
\psfrag{R3'}{$R_3'$}

     \includegraphics[width=15cm]{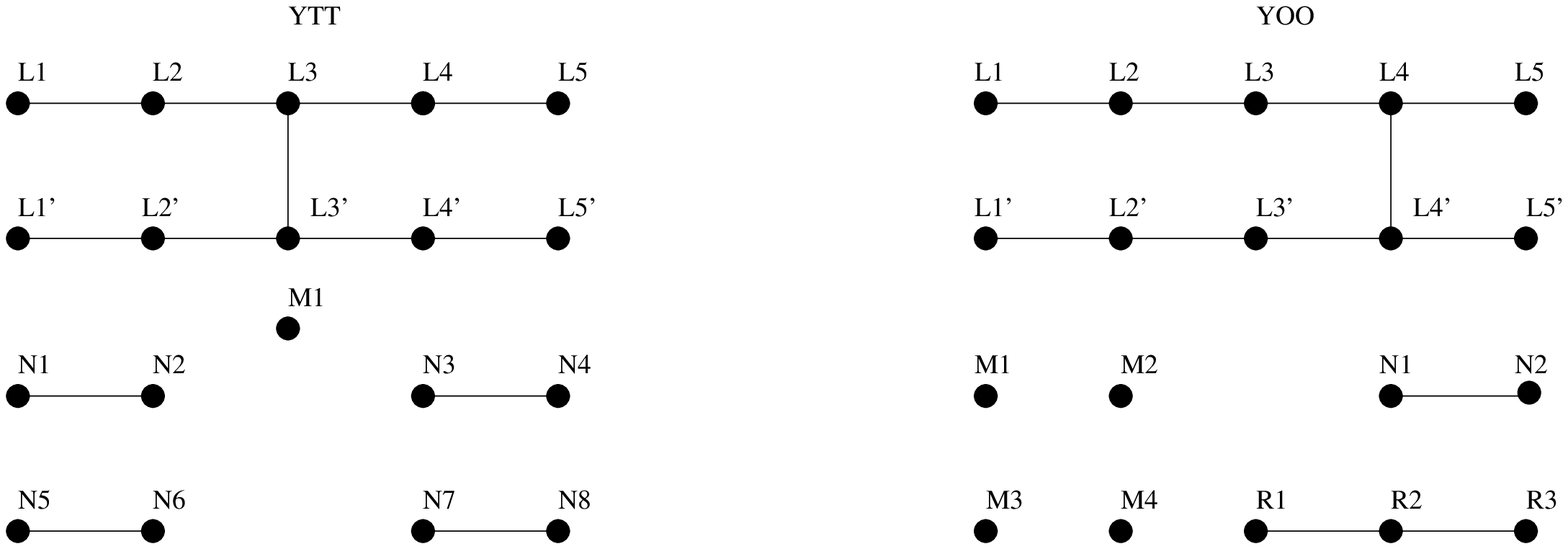}
   \end{psfrags}
   \vspace*{2mm}
\end{center}


\vspace*{3.0cm}





   \begin{center}
   \begin{psfrags}
     \psfrag{TTV}{$Y_{\lambda,T\times V}$}
     \psfrag{YOT}{$Y_{\lambda,O\times T}$}
     \psfrag{L1}{$L_1$}
     \psfrag{L2}{$L_2$}
     \psfrag{L3}{$L_3$}
     \psfrag{L4}{$L_4$}
     \psfrag{L5}{$L_5$}
\psfrag{L1'}{$L_1'$}
     \psfrag{L2'}{$L_2'$}
     \psfrag{L3'}{$L_3'$}
     \psfrag{L4'}{$L_4'$}
     \psfrag{L5'}{$L_5'$}
\psfrag{L1''}{$L_1''$}
     \psfrag{L2''}{$L_2''$}
     \psfrag{L3''}{$L_3''$}
     \psfrag{L4''}{$L_4''$}
     \psfrag{L5''}{$L_5''$}
\psfrag{M1}{$M_1$}
\psfrag{M2}{$M_2$}
\psfrag{M3}{$M_3$}
 \psfrag{L0}{$L_0$}
\psfrag{N1}{$N_1$}
\psfrag{N2}{$N_2$}
\psfrag{N3}{$N_3$}
\psfrag{N4}{$N_4$}
\psfrag{N5}{$N_5$}
\psfrag{N6}{$N_6$}
\psfrag{N7}{$N_7$}
\psfrag{N8}{$N_8$}
\psfrag{N9}{$N_9$}
\psfrag{N10}{$N_{10}$}
\psfrag{N11}{$N_{11}$}
\psfrag{N12}{$N_{12}$}
\psfrag{R2}{$R_{2}$}
     \includegraphics[width=15cm]{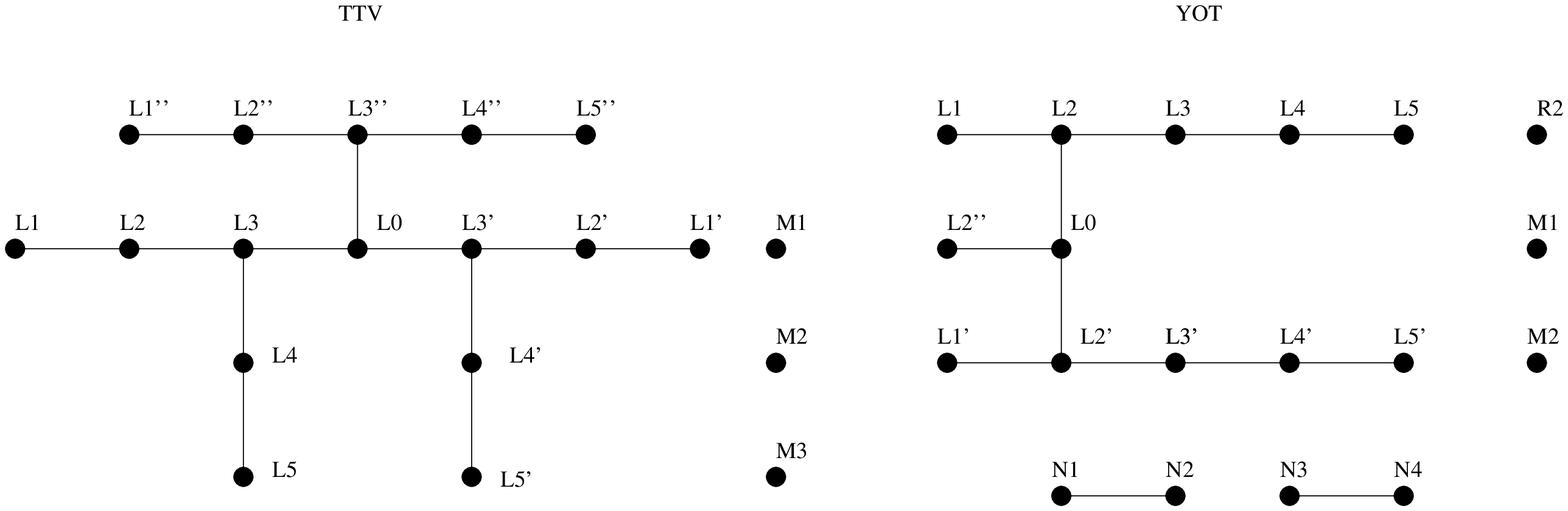}
   \end{psfrags}
   \vspace*{2mm}
\end{center}


\newpage




   \begin{center}
   \begin{psfrags}
\psfrag{61=64}{$Y^{(6,1)}_{\lambda,T\times V}(Y^{(6,4)}_{\lambda,T\times V})$}
\psfrag{62=63}{$Y^{(6,2)}_{\lambda,T\times V}(Y^{(6,3)}_{\lambda,T\times V})$}
\psfrag{81}{$Y^{(8,1)}_{\lambda,O\times T}$}
 \psfrag{82}{$Y^{(8,2)}_{\lambda,O\times T}$}
  \psfrag{83}{$Y^{(8,3)}_{\lambda,O\times T}$}
  \psfrag{84}{$Y^{(8,4)}_{\lambda,O\times T}$}

     \psfrag{L1}{$L_1$}
     \psfrag{L2}{$L_2$}
     \psfrag{L3}{$L_3$}
     \psfrag{L4}{$L_4$}
     \psfrag{L5}{$L_5$}
\psfrag{L1'}{$L_1'$}
     \psfrag{L2'}{$L_2'$}
     \psfrag{L3'}{$L_3'$}
     \psfrag{L4'}{$L_4'$}
     \psfrag{L5'}{$L_5'$}
\psfrag{L1''}{$L_1''$}
     \psfrag{L2''}{$L_2''$}
     \psfrag{L3''}{$L_3''$}
     \psfrag{L4''}{$L_4''$}
     \psfrag{L5''}{$L_5''$}
\psfrag{M1}{$M_1$}
\psfrag{M2}{$M_2$}
\psfrag{M3}{$M_3$}
 \psfrag{C}{$C$}
 \psfrag{L0}{$L_0$}
\psfrag{N1}{$N_1$}
\psfrag{N2}{$N_2$}
\psfrag{N3}{$N_3$}
\psfrag{N4}{$N_4$}
\psfrag{N5}{$N_5$}
\psfrag{N6}{$N_6$}
\psfrag{N7}{$N_7$}
\psfrag{N8}{$N_8$}
\psfrag{N9}{$N_9$}
\psfrag{N10}{$N_{10}$}
\psfrag{N11}{$N_{11}$}
\psfrag{N12}{$N_{12}$}
\psfrag{R2}{$R_{2}$}
     \includegraphics[width=15cm]{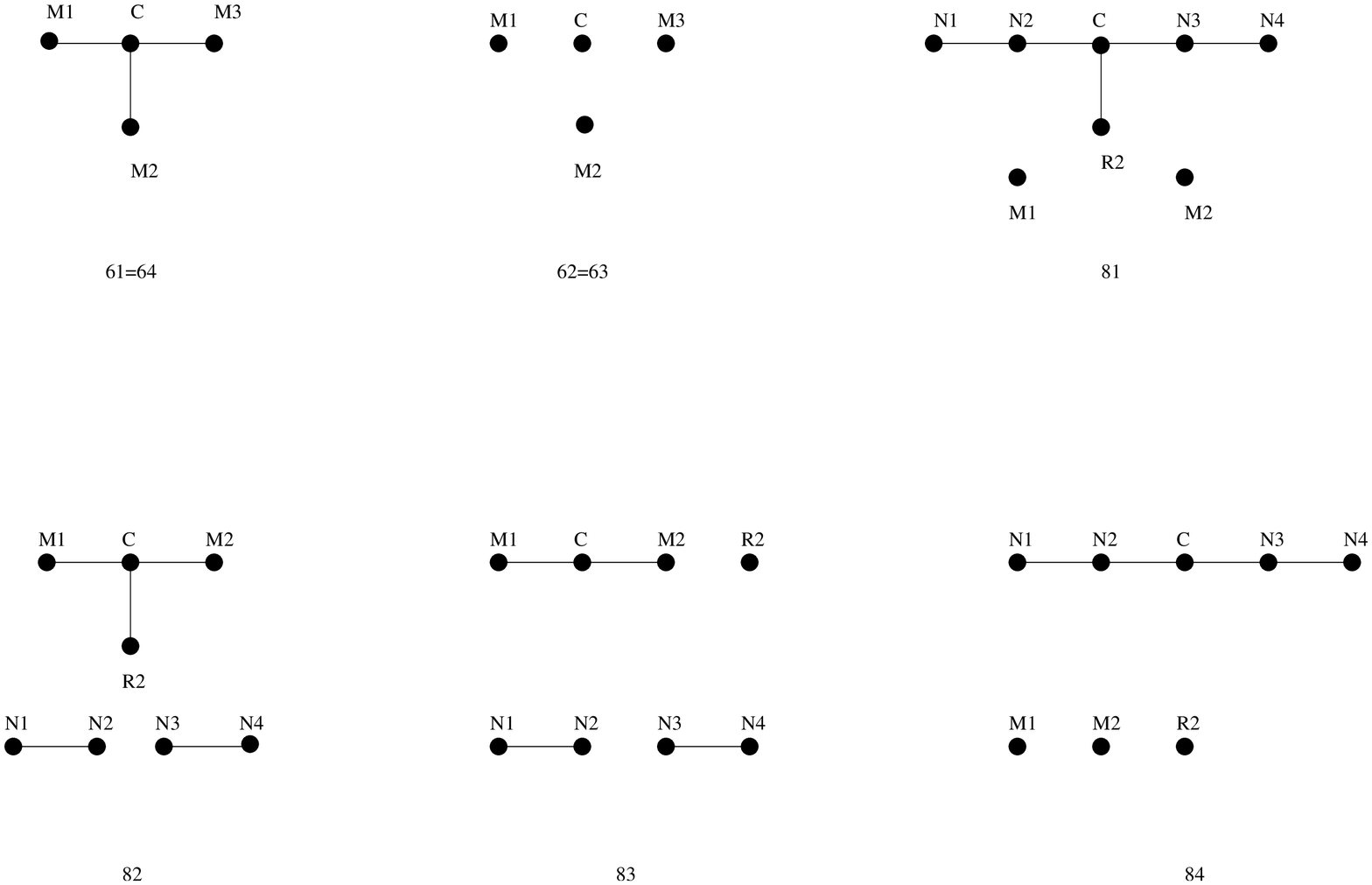}
   \end{psfrags}
   \vspace*{2mm}
\end{center}



\vspace*{0.2cm}




   \begin{center}
   \begin{psfrags}
\psfrag{61}{$Y^{(6,1)}_{\lambda,(TT)'}$}
\psfrag{62}{$Y^{(6,2)}_{\lambda,(TT)'}$}
\psfrag{81}{$Y^{(8,1)}_{\lambda,(OO)''}$}
 \psfrag{84}{$Y^{(8,4)}_{\lambda,(OO)''}$}
  \psfrag{811}{$Y^{(8,1)}_{\lambda,T\times T}$}
  \psfrag{844}{$Y^{(8,4)}_{\lambda,T\times T}$}

     \psfrag{L1}{$L_1$}
     \psfrag{L2}{$L_2$}
     \psfrag{L3}{$L_3$}
     \psfrag{L4}{$L_4$}
     \psfrag{L5}{$L_5$}
\psfrag{L1'}{$L_1'$}
     \psfrag{L2'}{$L_2'$}
     \psfrag{L3'}{$L_3'$}
     \psfrag{L4'}{$L_4'$}
     \psfrag{L5'}{$L_5'$}
\psfrag{L1''}{$L_1''$}
     \psfrag{L2''}{$L_2''$}
     \psfrag{L3''}{$L_3''$}
     \psfrag{L4''}{$L_4''$}
     \psfrag{L5''}{$L_5''$}
\psfrag{M1}{$M_1$}
\psfrag{M2}{$M_2$}
\psfrag{M3}{$M_3$}
\psfrag{M4}{$M_4$}
 \psfrag{C}{$C$}
 \psfrag{C1}{$C_1$}
 \psfrag{C2}{$C_2$}
 \psfrag{C3}{$C_3$}
 \psfrag{L0}{$L_0$}
\psfrag{N1}{$N_1$}
\psfrag{N2}{$N_2$}
\psfrag{N3}{$N_3$}
\psfrag{N4}{$N_4$}
\psfrag{N5}{$N_5$}
\psfrag{N6}{$N_6$}
\psfrag{N7}{$N_7$}
\psfrag{N8}{$N_8$}
\psfrag{N9}{$N_9$}
\psfrag{N10}{$N_{10}$}
\psfrag{N11}{$N_{11}$}
\psfrag{N12}{$N_{12}$}
\psfrag{R1}{$R_{1}$}
\psfrag{R2}{$R_{2}$}
\psfrag{R3}{$R_{3}$}
\psfrag{R1'}{$R_{1}'$}
\psfrag{R2'}{$R_{2}'$}
\psfrag{R3'}{$R_{3}'$}

     \includegraphics[width=15cm]{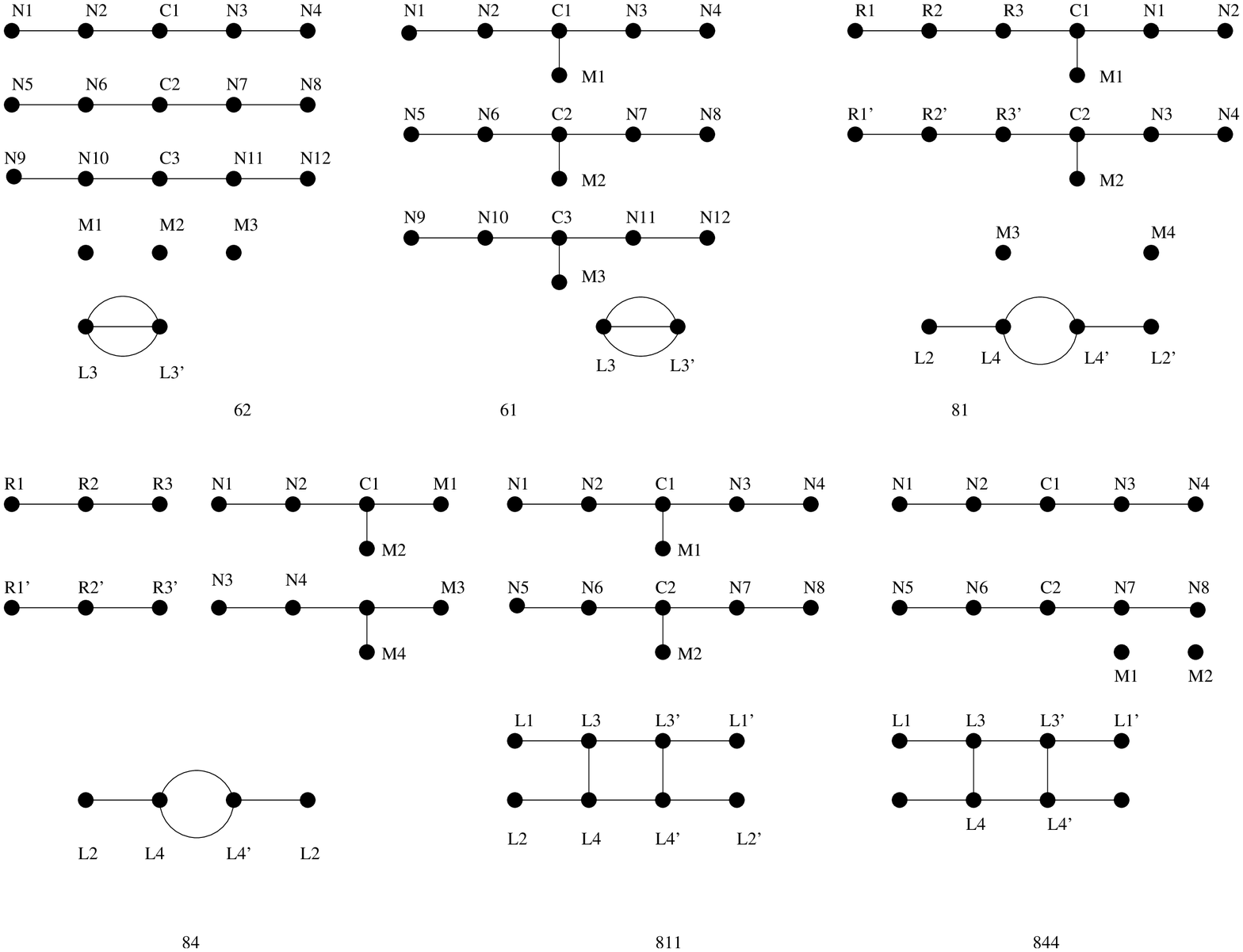}
   \end{psfrags}
   \vspace*{2mm}
\end{center}



\smallskip


\addcontentsline{toc}{section}{  \hspace{0.5ex} References}

\end{document}